\documentclass{amsart}

\usepackage{graphicx}

\newtheorem{theorem}{Theorem}[section]
\newtheorem{lemma}[theorem]{Lemma}
\newtheorem{proposition}[theorem]{Proposition}
\theoremstyle{definition}
\newtheorem{definition}[theorem]{Definition}
\newtheorem{example}[theorem]{Example}

\newtheorem{problem}[theorem]{Problem}
\theoremstyle{remark}
\newtheorem{remark}[theorem]{Remark}
\newtheorem{acknowledgement}[theorem]{Acknowledgement}
\numberwithin{equation}{section}

\begin{document}

\title{Internal Polya Inequality for $\mathbb{C}$-convex domains in $\mathbb{C}^{n}$}

\author{OZAN G\"UNY\"UZ}
\address{Department of Mathematics, Sabanc{\i} University, 34956 Tuzla/Istanbul, Turkey}
\email{ozangunyuz@sabanciuniv.edu}

\author{VYACHESLAV ZAKHARYUTA}
\address{Department of Mathematics, Sabanc{\i} University, 34956 Tuzla/Istanbul, Turkey}

\email{zaha@sabanciuniv.edu}


\subjclass[2010]{Primary 54C40, 14E20; Secondary 46E25, 20C20}

\begin{abstract}
Let $K\subset \mathbb{C}$ be a polynomially convex compact set, $f$ be a
function analytic in a domain $\overline{\mathbb{C}}\backslash K$ with
Taylor expansion $f( z) =\sum_{k=0}^{\infty }\frac{a_{k}}{z^{k+1}}
$ at $\infty $, and $H_{i}( f) :=\det ( a_{k+l})
_{k,l=0}^{i}$ the related Hankel determinants. The classical Polya theorem \cite%
{P} says that
\[
\limsup_{i\rightarrow \infty }\vert H_{i}( f) \vert
^{1/i^{2}}\leq d( K) ,
\]%
where $d( K) $ is the transfinite diameter of $K$. The main result of this paper is a multivariate analog of Polya's inequality for a weighted Hankel-type determinant constructed from the Taylor series of a function analytic on a $\mathbb{C}$-convex (=strictly linearly convex) domain in $\mathbb{C}^{n}$.
\end{abstract}

\maketitle

\section*{Introduction}

Suppose $f$ is a function analytic in a domain $\overline{\mathbb{C}}\backslash K$, where $K\subset \mathbb{C}$ is a polynomially convex
compact set, with Taylor expansion $f( z) =\sum_{k=0}^{\infty }%
\frac{a_{k}}{z^{k+1}}$ at $\infty $, and $H_{i}( f) :=\det (
a_{k+l}) _{k,l=0}^{i}$ are the related Hankel determinants. The classical
Polya theorem \cite{P} says that
\[
D( f) :=\limsup_{i\rightarrow \infty }\vert H_{i}(
f)\vert ^{1/i^{2}}\leq d( K) ,
\]%
where $d( K) $ is the transfinite diameter of $K$ introduced by
Fekete in \cite{F}.

The present paper studies some multivariate analogs of this result. Leja
\cite{Lj} introduced the notion of tranfinite diameter for a compact set $K\subset \mathbb{C}^{n}$. In order to represent his definition\ we need some
notation. Let
$$ e_{i}( z) :=z^{k( i) }=z_{1}^{k_{1}( i)
}\ldots z_{\nu }^{k_{\nu }( i) }\ldots z_{n}^{k_{n}(
i) },\ i\in \mathbb{N} $$
be all monomials in $\mathbb{C}^{n}$, ordered so that the degrees $\sigma(i) :=k_{1}(i)+\ldots+k_{n}(i)$ be non-decreasing. The number of all monomials of degree $\leq \sigma$ is  $m_{\sigma}$
and $l_{\sigma}:=\sum_{q=1}^{\sigma}q( m_{q}-m_{q-1})$. The \textit{transfinite diameter} of $K$ is the number%
\begin{equation}
d( K) :=\limsup_{\sigma\rightarrow \infty }( \max \{
\vert \det ( e_{i}( \zeta _{j}) )
_{i,j=1}^{m_{\sigma}}\vert :( \zeta _{j}) \in K^{m_{\sigma}}\}
) ^{1/l_{\sigma}}.  \label{dk}
\end{equation}

In particular, this definition coincides with that of Fekete for $n=1$; it
was proved (for $n=1$ in \cite{F} and for $n\geq 2$ in \cite{Z1}) that the
usual limit exists in (\ref{dk}).

A direct analog of Polya's inequality has no sense, since there are no
bounded analytic functions except identical constants in the complement of a
polynomially convex compact set in $\overline{\mathbb{C}^{n}}$. Schiffer and
Siciak \cite{SS} considered an analog of Polya's inequality in terms of $%
d( K) $ for a Cartesian product of plane polynomially convex
compact sets $K=K_{1}\times \ldots \times K_{n}$ and functions analytic in
the conjugate polycylindrical domain $D=( \overline{\mathbb{C}}\backslash K_{1}) \times \ldots \times ( \overline{\mathbb{C}}\backslash K_{n}) $, which plays the role of an "exterior" for $K$.
Zakharyuta \cite{Z1} obtained a general multivariate version of Polya's
inequality with $d( K) $ for a polynomially convex compact set $%
K\subset \mathbb{C}^{n}$ and an \textit{analytic functional} $F$ in $\mathbb{%
C}^{n}$ which can be extended as a linear continuous functional onto the
space $A( K) \supset A( \mathbb{C}^{n}) $ of all germs
of analytic functions on $K$, endowed with the natural inductive topology( see Section 2 below).
Therewith the related Hankel determinants are built from the Taylor
coefficients of the analytic germ $\varphi $ at the point $( \infty
,\ldots, \infty ) \in \overline{\mathbb{C}}^{n}$, generated by the
analytic functional $F$ (for more detail see below, Section 2).

We consider in Section 3 multivariate internal analogs of Polya's inequality
for a strictly linearly convex domain $D\subset \mathbb{C}^{n}$ and
Hankel-type determinants, constructed from the Taylor coefficients of a
function $f\in A( D) $ at a given point $a\in D$; therewith these determinants include special weights, generated by $s$-indicatrices (see Section 3) of the sequence of analytic functionals biorthogonal to the system of monomials in $\mathbb{C}^{n},\,\,s=1,\ldots,n$. The
main result (Theorem \ref{ipt}) is obtained by a reduction to the aforementioned result from \cite{Z1}. It is based on S. Znamenskii's characterization of strict linear convexity in terms of $s$-indicatrices
(see, Theorem \ref{tzn} below).

V. Sheinov has claimed an analogous result for $s=n$, under more
restrictive assumptions about $D$. We analyse his results in Section 3 and
find out that the related statements (\cite{Sh1},Theorem 2 and \cite{Sh2},
Theorem A) are true only if $a=0$, otherwise they are wrong, in general,
even in the case $n=1$.

In section 3, we provide a counterexample to Sheinov's claim and give several
examples of internal Polya inequality considering different sets in $\mathbb{C}^n.$

The internal transfinite diameter $d(0,\partial D)$ (see, Definition 4.1) is considered in
Section 4 in the context of Theorem 2.4 and internal Polya inequality for nonweighted
Hankel determinants (see, Problem 4.5 and a conjecture after it).

\section{Polya's Theorem in $\mathbb{C}^{n}$}
By $A( D) $ we mean the space of all analytic functions in a
domain $D\subset \mathbb{C}^{n}$ with usual locally convex topology of
locally uniform convergence in $D$. For a compact set $K\subset \mathbb{C}%
^{n}$ we denote by $A( K) $ the locally convex space of all germs
of analytic functions on $K$, endowed with the standard inductive topology.
Elements of the dual space $A( \mathbb{C}^{n}) ^{\ast }$are
called \textit{analytic functionals in }$\mathbb{C}^{n}$ (see, e.g., H\"{o}%
rmander \cite{Hrm}, \ Section 4.5). By $A_{o}( \{ \infty
^{n}\} ) $\textit{\ }we denote the space of all germs of
analytic functions $\varphi $ at the point $\infty ^{n}:=( \infty
,\ldots ,\infty ) \in \overline{\mathbb{C}^{n}}:=\overline{\mathbb{C}}%
\times ,\ldots ,\times \overline{\mathbb{C}}$, having an expansion :
\begin{equation}
\varphi ( z) =\sum_{k\in \mathbb{Z}_{+}^{n}}\frac{a_{k}(
\varphi ) }{z^{k+I}},\,\,\,I:=( 1,\ldots ,1)  \label{expfi}
\end{equation}%
converging uniformly in a neighborhood $\{ z=( z_{\nu }) \in
\overline{\mathbb{C}^{n}}:\vert z_{\nu }\vert \geq r\} $ of
$\infty ^{n}$ with some $r=r( \varphi ) $. The following fact is
well-known (see, e.g., \cite{Dn}).

\begin{lemma}
There is an isomorphism $T:A( \mathbb{C}^{n}) ^{\ast }\rightarrow
A_{o}( \{ \infty ^{n}\} ) $ such that for every $F\in
A( \mathbb{C}^{n}) ^{\ast }$ and $\varphi =T( F) $
there exists $r\,=r( F) >0$ so that the formula %
\begin{equation}
F( f) =[ f,\varphi ] :=( \frac{1}{2\pi i})
^{n}\int_{\mathbb{T}_{R}^{n}}f( \zeta ) \varphi ( \zeta
) d\zeta ,  \label{ff}
\end{equation}%
holds for each $f\in A( \mathbb{C}^{n}) $ and any $R>r$; here%
\begin{equation}
\mathbb{T}_{R}^{n}:=\{ z=( z_{\nu }) :\vert z_{\nu
}\vert =R,\ \nu =1,\ldots ,n\} .  \label{trn}
\end{equation}
\end{lemma}

This lemma provides, for every analytic functional $F$, a related sequence
of multivariate Hankel-like determinants:%
\begin{equation}
H_{i}=H_{i}( \varphi ) :=\det(a_{k( \alpha)+k(\beta) })_{\alpha ,\beta =1}^{i},\ i\in \mathbb{N},
\label{hif}
\end{equation}%
with
\[
a_{k( \alpha ) }=F( e_{\alpha }) =[ e_{\alpha
},\varphi ] ,\ \varphi =T( F) ,\ \alpha \in \mathbb{N}.
\]%
The general multivariate analog of Polya's Theorem is the following

\begin{theorem}
\label{mpt}(\cite{Z1}) Suppose $K$ is a polynomially convex compact set in $%
\mathbb{C}^{n}$, $F$ is an analytic functional which has a continuous
extension onto $A( K) $ and $\varphi =T( F) $ is the
corresponding analytic germ at $\infty ^{n}$. Then the determinants (\ref%
{hif}) satisfy the following asymptotic inequality:%
\[
D( \varphi ) :=\limsup\limits_{i\rightarrow \infty }\vert
H_{i}( \varphi ) \vert ^{\frac{1}{2l_{\sigma( i) }}%
}\leq d( K) .
\]%
\qquad \qquad
\end{theorem}

It has been proved in \cite{Z1} a bit weaker statement with the outer
transfinite diameter ; $\widehat{d}( K) $ instead of $d(
K) $; it was proved later that $\widehat{d}( K) =d(
K) $ (\cite{Znv1,Znv2,Lv}).

\section{Internal Polya theorem}
Notice that for $n=1$ the following internal version of Polya Theorem can be
easily derived from the classical result of Polya.

\begin{proposition}
\label{1ipineq}Let $D\subset \mathbb{C}$ be a domain, $a\in D$. Let $f$ be a
function analytic in $D$ with its Taylor expansion at the point $a$:%
\[
f( z) =\sum_{k=0}^{\infty }c_{k}\ ( z-a) ^{k}
\]%
and $H_{i}( f,a) :=\det ( c_{k+l}) _{k,l=0}^{i}$ be
the related Hankel determinants. Then
\begin{equation}
  \limsup\limits_{i\rightarrow \infty }\vert H_{i}( f,a)
\vert ^{\frac{1}{i^{2}}}\leq d( K_{a}) ,
\label{ipo}\end{equation}
where $K_{a}=\{\frac{1}{z-a}:z\in \overline{\mathbb{C}}
\backslash D\}$.
\end{proposition}

We are going to prove a multivariable analog of this fact for the class of
strictly linearly convex domains in $\mathbb{C}^{n}$ (known also as $\mathbb{%
C}$-convex). First we give necessary definitions. For more detail, see \cite{APS} or \cite{Kis}.

In what follows, $E$ denotes a domain or a compact set in $\mathbb{C}^{n}$.
The set
\[
E^{\ast }:=\{ w=( w_{\nu }) \in \mathbb{C}^{n}:w\cdot z\neq 1%
\,\,\, for\,\,\, all\,\,\, z=( z_{\nu }) \in E\} ,
\]%
where $w\cdot z:=\sum_{\nu =1}^{n}w_{\nu }z_{\nu },$ is called the \textit{%
conjugate to }$E$.

A set $E$ is called \textit{linearly convex }if $E^{\ast \ast }:=(
E^{\ast }) ^{\ast }=E$. This is equivalent to saying that $E$ is linearly convex if its complement is a union of complex hyperplanes. If $E$ is linearly convex then the functions

\[
u_{w}^{( s) }( z) :=( 1-wz) ^{-s},\ z\in
E;\ w\in E^{\ast },\ s=1,\ldots ,n
\]%
belong to the space $A( E) $. Given a \ functional $F\in $%
\bigskip $A( E) ^{\ast }$, the function
\[
\varphi _{s}( w) =F( u_{w}^{( s) })
\]%
belongs to the space $A( E^{\ast }) $ and is called the $s$%
-indicatrix of the functional $F$. This correspondence defines a linear
operator
\[
\mathcal{F}_{s}:A( E) ^{\ast }\rightarrow A( E^{\ast
})
\]%
for each $s=1,\ldots ,n$. The operator $\mathcal{F}:=\mathcal{F}_{1}$ is
known as the \textit{Fantappie transform}.

\begin{definition}
A linearly convex set $E\subset \mathbb{C}^{n}$ is said to be strictly
linearly convex, if the \textit{Fantappie transform }$\mathcal{F}:A(
E) ^{\ast }\rightarrow A( E^{\ast }) $ is an isomorphism; $%
E $ is called $\mathbb{C}$-convex, if for each complex line $l$ the set $%
l\cap E$ is connected and simply connected.
\end{definition}

\begin{theorem}
\label{tzn}(S. Znamenskii \cite{Zn}, see also \cite{APS}) Given a linearly
convex set $E\subset \mathbb{C}^{n}$, the following statements are
equivalent: $( i) $ $E$ is strictly linearly convex, $(
ii) $ $E$ is $\mathbb{C}$-convex, $( iii) $ The operator $%
\mathcal{F}_{s}$ is an isomorphism for each $s=1,\ldots ,n,$ $(
iv) $ The operator $\mathcal{F}_{s}$ is an isomorphism for some $%
s=1,\ldots ,n$.
\end{theorem}

Let%
\[
h_{k}( w) =w^{k};\ \ h_{k}^{\prime }\in A( \mathbb{C}%
^{n}) ^{\ast },\ h_{k}^{\prime }( h) :=\frac{h^{(
k) }( 0) }{k!}
\]%
be the system of monomials in $\mathbb{C}_{w}^{n}$ and its biorthogonal
system, respectively. Then, by a simple calculation%
\begin{equation}
\mathcal{F}_{s}( h_{k}^{\prime }) ( z) =\frac{(
u_{z}^{( s) }( 0) ) ^{( k) }}{k!}%
=\lambda _{k}^{( s) }\ z^{k},\ k\in \mathbb{Z}_{+}^{n},
\label{diag}
\end{equation}%
where
\begin{equation}
\lambda _{k}^{( s) }\ =\frac{(\vert k \vert +s-1) !}{k!(s-1)!}  \label{lk}
\end{equation}

\begin{theorem}
\label{ipt}Let $D$ be a strictly linearly convex domain in $\mathbb{C}^{n}$,
$a\in D$. Let $f\in A( D) $,
\begin{equation}
f( z) =\sum_{k\in \mathbb{Z}_{+}^{n}}c_{k}\ ( z-a) ^{k}
\label{expf}
\end{equation}%
be its Taylor expansion at $z=a$, and
\[
H_{i}^{( s) }( f,a) =\det ( \frac{c_{k(
\alpha ) +k( \beta ) }}{\lambda _{k( \alpha )
+k( \beta ) }^{( s) }}) _{\alpha ,\beta
=1}^{i},\ i\in \mathbb{N},
\]%
be weighted Hankel-type determinants. Then%
\begin{equation}
\limsup\limits_{i\rightarrow \infty }\vert H_{i}^{( s)
}( f,a) \vert ^{1/2l_{\sigma( i) }}\leq d(
( D-a) ^{\ast })   \label{wpineq}
\end{equation}
for each $s=1,\,2,\,\ldots,\,n$.
\end{theorem}

\begin{proof}
Since the operator $T:A( D) \rightarrow A( D-a) $,
defined by%
\[
( Tf) ( w) :=f( w-a) ,
\]%
is an isomorphism, it is sufficient to consider only the case $a=0$. By the
S. Znamenskii result (Theorem \ref{tzn} above) each isomorphism $\mathcal{F}%
_{s}:A( D^{\ast }) ^{\ast }\rightarrow A( D^{\ast \ast
}) =A( D) $ is a continuous linear extension of the
diagonal mapping (\ref{diag}); it maps the functional $h^{\prime }\in
A( D^{\ast }) ^{\ast }$ with the formal expansion $h^{\prime
}\sim \sum_{k}h^{\prime }( h_{k}) h_{k}^{\prime }$ onto the
function $f\in A( D) $ with the Taylor expansion
\[
f( z) =\sum_{k\in \mathbb{Z}_{+}^{n}}h^{\prime }(
h_{k}) \ \lambda _{k}^{( s) }z^{k},\ s=1,\ldots ,n,
\]%
converging in some open neighborhood of zero. Conversely, if $f\in A(
D) $ has the Taylor expansion (\ref{expf}) with $a=0$, then the
analytic functional
\[
h^{\prime }=\mathcal{F}_{s}^{-1}( f) =\sum_{k\in \mathbb{Z}%
_{+}^{n}}\frac{c_{k}}{\lambda _{k}^{( s) }}\ h_{k}^{\prime }\
\]%
is extendible continuously onto the space \,$A( D^{\ast }) $\, and the
corresponding analytic germ $\varphi $ has the expansion (\ref{expfi}) at $%
\infty ^{n}$ with $a_{k}( \varphi ) =\frac{c_{k}}{\lambda
_{k}^{( s) }}$. Hence $G_{i}^{( s) }( f,0)
=H_{i}( \varphi ) $ and, by Theorem \ref{mpt}, we have (\ref%
{wpineq}) with $a=0$.
\end{proof}

\begin{remark}
Proposition \ref{1ipineq} is a particular case of Theorem \ref{ipt}, since
every domain $D\subset \mathbb{C}$ is strictly linearly convex,%
\[
( D-a) ^{\ast }=K_{a}=\{ \frac{1}{z-a}:z\in \overline{%
\mathbb{C}}\backslash D\},\,\,\ a\in D,
\]%
and the weight $\lambda _{k}^{( 1) }$, for $n=1,$ equals $1$
identically.

In contrast to Theorem \ref{mpt} and Proposition \ref{1ipineq}, Theorem \ref{ipt} is proved under a quite strict assumption about the domain $D$, therewith the Hankel-like determinants contain special weights. Those weights are essential and cannot be removed. The question on the internal Polya inequality for non-weighted Hankel determinants is still open (see the discussion at the end of the paper).
\label{remark} \end{remark}

\section{Examples}

In this section, we will present some examples concerning Theorem \ref{ipt}.

\begin{example}
Let $D=\{ z\in \mathbb{C}:\vert z\vert <1\} ,$ $a\in
D,\ a\neq 0$. It is easy to see that $( D-a) ^{\ast }=K_{a}$ is a
disc of radius $\frac{1}{1-\vert a\vert ^{2}}$, hence $%
d( K_{a}) =\frac{1}{1-\vert a\vert ^{2}}$. By
the classical Polya inequality and Goluzin's result on sharpness of Polya
inequality (\cite{G}, XI.5, Theorem 2), there exists a function $g_{a}\in
A( \overline{\mathbb{C}}\backslash K_{a}) $, having the
Taylor expansion
\[
g_{a}( z) =\sum_{k=0}^{\infty }\frac{c_{k}^{( a) }}{%
z^{k+1}},
\]%
and such that its Hankel determinants $H_{i}( g_{a}) :=\det
( c_{k+l}^{( a) }) _{k,l=0}^{i}$ satisfy the
asymptotic inequality
\begin{equation}
1<\limsup\limits_{i\rightarrow \infty }\vert
H_{i}( g_{a}) \vert ^{1/i^{2}}\leq d( K_{a}) .
\label{1<d}
\end{equation}%
Since the operator $S:A( \overline{\mathbb{C}}\backslash
K_{a}) \rightarrow A( D) $, defined by the formula:
\[
S( g) ( \zeta ) :=\frac{1}{\zeta -a}g( \frac{1}{%
\zeta -a}) ,
\]%
is an isomorphism, the function $f_{a}=S( g_{a}) $ belongs to the
space $A( D) $ and has its Taylor expansion at $\zeta =a$ with
the same coefficients as $g_{a}$:
\[
f_{a}( \zeta ) =\sum_{k=0}^{\infty }c_{k}^{( a)
}( \zeta -a) ^{k}.
\]%
Therefore $H_{i}( f_{a},a) =H_{i} g_{a}) $ and (\ref%
{1<d}) implies
\[
d( D^{\ast }) =1<\limsup\limits_{i\rightarrow \infty }\vert
H_{i}( f_{a},a) \vert ^{\frac{1}{i^{2}}}\leq d(
K_{a}).
\]%

\end{example}

Sheinov stated in \cite{Sh1, Sh2} (for $s=n$ and under some additional,
quite restrictive assumptions), that the asymptotic estimate (\ref{wpineq})
holds with $d( D^{\ast }) $ instead of $d( ( D-a)
^{\ast }) $. As the  above example shows, his statement is wrong
if $a\neq 0$ even for the plane unit disc (we are not able to comment on his proof,
since it is inaccessible).

\begin{example}

Let \begin{equation}
D=\{ w=( w_{\nu }) :\sum_{\nu =1}^{n}\vert w_{\nu}\vert R_{\nu }< 1\}.
                                                                  \label{dual}
                                                                  \end{equation} be the $n$-dimensional open  hypercone. Then $D^{*}$ is the closed  polydisc with polyradii $R_{\nu},\,\,\nu=1,\ldots,n.$ The transfinite diameter of this set (see for example, \cite{SS}) is $$d(D^{*})=\sqrt[n]{R_1.\ldots.R_n}.$$ By Theorem \ref{ipt}, we have $$\limsup\limits_{i\rightarrow \infty }\vert H_{i}^{( s)
}( f,0) \vert ^{1/2l_{\sigma( i) }}\leq \sqrt[n]{R_{1}.\ldots.R_{n}}. $$
\end{example}
\begin{example}
Let now $D$ be the open polydisc with polyradii $R_{\nu},\,\,\nu=1,\ldots,n.$ Since $D^{*}$ is, in this case,  the $n$-dimensional closed hypercone, it follows from the result of M. Jedrzejowski (\cite{Jed}) and Theorem \ref{ipt}  $$\limsup\limits_{i\rightarrow \infty }\vert H_{i}^{( s)
}( f,0) \vert ^{1/2l_{\sigma( i) }}\leq \frac{\exp ( 1+\sum_{i=1}^{n}(
-1) ^{i}C_{n}^{i}\frac{1}{i}) }{ \sqrt[n]{R_{1}\cdots
R_{n}}}.$$ \label{expp}

\end{example}

We will give an alternate proof of the transfinite diameter of closed hypercone. For this, we need to give some information. Given that $K$ is an $n$-circular complete compact set in $\mathbb{C}^n,$ let $h_{K}:\Sigma \rightarrow \mathbb{R}$ be its characteristic function defined by  \begin{equation} h_{K}( \theta ):=\max \{ \sum_{\nu =1}^{n}\theta
_{\nu }\ln \vert w_{\nu }\vert :w=( w_{\nu }) \in
K\}. \label{max} \end{equation} where $\lambda$ is the standard $(n-1)$-dimensional simplex
\begin{equation}
\Sigma :=\{ \theta =(\theta_{\nu }) \in \mathbb{R}%
^{n}:\theta _{\nu }\geq 0,\ \nu =1,\ldots ,n;\ \sum_{\nu =1}^{n}\theta _{\nu
}=1\}.  \label{sgm}
\end{equation}
By [\cite{Z1}, section 7], we have the following relation \begin{equation} d(K)=\exp{(\frac{1}{\lambda{}(\Sigma)}\int_{\Sigma} (h_{K}(\theta)) d\lambda ( \theta )}) .
 \end{equation} where $\lambda$ is the Lebesgue surface area measure on $\Sigma.$
 \begin{lemma}  Let $K$ be the $n$-dimensional closed hypercone $\{ w=( w_{\nu }) :\sum_{\nu =1}^{n}\vert w_{\nu}\vert R_{\nu }\leq 1\}.$ Then we have
$$d(K) =\frac{\exp ( 1+\sum_{i=1}^{n}(
-1) ^{i}C_{n}^{i}\frac{1}{i}) }{( R_{1}\cdots
R_{n}) ^{1/n}}.$$
     \label{hyperlem} \end{lemma}
\begin{proof}
  Initially, let us determine the characteristic function of $K$\,\,explicitly. Since the maximum is attained on the set $$ \{ w=( w_{\nu }) :\sum_{\nu =1}^{n}\vert w_{\nu
}\vert R_{\nu }= 1\}, $$
writing $\xi_{\nu}:=\ln{\vert {w_{\nu}} \vert}$, (\ref{max}) can be written as follows $$h_{K}(\theta)=\max \{ \sum_{\nu =1}^{n}\theta _{\nu }\xi _{\nu }:\sum_{\nu
=1}^{n}R_{\nu }\exp \xi _{\nu }=1\}.$$
By an application of Lagrange multipliers to the summation $\sum_{\nu=1}^{n}{\theta_{\nu}\xi_{\nu}}$ subject to the constraint $\sum_{\nu=1}^{n}{R_{\nu}\exp{\xi_{\nu}}}=1$, we have $$h_{K}(\theta)=\sum_{\nu=1}^{n}{\theta_{\nu}\ln{\frac{\theta_{\nu}}{R_{\nu}}}}.$$ Therefore, \begin{eqnarray}
d(K) &=&\exp \frac{1}{\lambda ( \Sigma ) }%
\int_{\Sigma }\sum_{\nu=1}^{n} \theta _{\nu }\ln \frac{\theta _{\nu }}{R_{\nu }}d\lambda
( \theta )  \nonumber \\
&=&\exp ( M_{1}+M_{2}) ,\label{plus}
\end{eqnarray}
where%
\[
M_{1}:=-\frac{1}{\lambda ( \Sigma ) }\int_{\Sigma }\sum_{\nu=1}^{n} \theta
_{\nu }\ln R_{\nu }d\lambda ( \theta ) =-\frac{1}{n}\sum_{\nu=1}^{n} \ln
R_{\nu }
\]%
and thus
\[
\exp M_{1}=\frac{1}{( R_{1}\cdots R_{n}) ^{1/n}}
\]%
is the transfinite diameter of the polydisc with polyradii $\frac{1}{R_{\nu }%
},\,\,\nu=1,\ldots,n$.

On the other hand,%
$$
M_{2}:=\frac{1}{\lambda ( \Sigma ) }\int_{\Sigma }\sum \theta
_{\nu }\ln \theta _{\nu }\ d\lambda ( \theta )
$$%
and thus $\exp M_{2}$ is the transfinite diameter of the unit hypercone%
$$
\{ z=( z_{\nu }) \in \mathbb{C}^{n}:\sum \vert z_{\nu
}\vert \leq 1\} .
$$%
By symmetry,%
$$M_{2}:=\frac{n}{\lambda ( \Sigma ) }\int_{\Sigma }\theta _{1}\ln
\theta _{1}\ d\lambda ( \theta )$$
We compute this surface integral as a multiple integral over the set%
$$
\Delta =\{ t=( t_{\nu }) \in \mathbb{R}_{+}^{n-1}:\sum
t_{\nu }\leq 1\}$$ by using the relations $d\lambda ( \theta ) =\sqrt{n}dt$ and $\lambda ( \Sigma
) =\frac{\sqrt{n}}{( n-1) !}$,\,where $dt$ is the Lebesgue measure on $\Delta.$

Then
\begin{eqnarray*}
M_{2} &=&n!\int_{\Delta }t_{1}\ln t_{1}\ dt  \nonumber \\
&=&n!\int_{0}^{1}t_{1}\ln t_{1}\ dt_{1}\int_{0}^{1-t_{1}}dt_{2}\cdots
\int_{0}^{1-t_{1}-\cdots -t_{n-2}}dt_{n-1} \\
&=&n( n-1) \int_{0}^{1}t_{1}\ln t_{1}\ ( 1-t_{1}) ^{n-2}dt\\
&=&1+\sum_{i=1}^{n}{(-1)^{i}C_{n}^{i}\frac{1}{i}}
\end{eqnarray*}
Finally taking the exponential and using the relation (\ref{plus}) will give the desired result.
\end{proof}

\section{Internal Transfinite Diameter}

First we recall the definition of the internal transfinite diameter,
following \cite{Zint}. Given a point $c\in \mathbb{C}^{n}$,
consider the system of analytic functionals $\{ e_{i ,a}^{\prime
}\} $, biorthogonal to the system of monomials $\{ (
z-c) ^{k(i)}\} $ and defined by%
\[
e_{i ,c}^{\prime }( f) =\frac{f^{k( i) }(
c) }{k( i) !},\ i\in \mathbb{N},\ f\in A( \{
c\} ) .
\]

\begin{definition}
\label{4.1}Let $D$ be a domain in $\mathbb{C}^{n}$, $c\in D.$ The
transfinite diameter of the boundary $\partial D$, viewed from the point $a$%
, is the quantity%
\[
d( c,\partial D) :=\limsup\limits_{i\rightarrow \infty }(
L_{i}) ^{\frac{1}{l_{\sigma ( i) }}},
\]%
where $L_{i}=\sup \{ \det ( e_{\alpha ,c}^{\prime }( f_{\beta
})) _{\alpha ,\beta =1}^{i}:f_{\beta }\in \mathbb{B}\} $
and $\mathbb{B}$ is the unit ball of the space $H^{\infty }( D) $ of bounded analytic functions.
\end{definition}

\begin{remark}
\label{4.2}It is proved in \cite{Zint} that in the one dimensional
case $d( c,\partial D) =d( K_{c}) $, where $K_{c}$ is
as in Remark \ref{remark}.
\end{remark}

\begin{example}
\label{4.3}Let $D$ be an open ball in $\mathbb{C}^{n}$ centered at $0$ with
radius $R$. It is known that its conjugated set $D^{\ast }$ is the closed
ball with radius $\frac{1}{R}$, hence $d( D^{\ast }) =\frac{1}{R}$. Due to \cite{Zint}, we have in this case also that $d(
0,\partial D) =\frac{1}{R}$.
\end{example}

Looking at these cases one can expect that Theorem \ref{ipt} is true with $d(
a,\partial D) $ instead of $d( ( D-a) ^{\ast }) $ but it is not the case, as it follows from the next.

\begin{example}
\label{4.4}Let $D$ be an open polydisc of polyradii $R_{1},\ldots ,R_{n}$,
centered at the origin. Then $D^{\ast }$ is the hypercone and due to Lemma \ref{hyperlem},
\[
d( D^{\ast }) =\frac{\exp ( 1+\sum_{i=1}^{n}( -1)
^{i}C_{n}^{i}\frac{1}{i}) }{\sqrt[n]{R_{1}\cdot \ldots \cdot R_{n}}},
\]%
while, by \cite{Zint},
\[
d( 0,\partial D) =\frac{1}{\sqrt[n]{R_{1}\cdot \ldots \cdot R_{n}}%
}.
\]%
Thus the diameter $d( c,\partial D) $ fails to be a proper
characteristic in the context of Theorem \ref{ipt}. However it is promising to deal
properly with the internal Polya inequality for non-weighted Hankel determinants,
as one can see below.
\end{example}

\begin{problem}
Let $D$ be a domain in $\mathbb{C}^{n}$, $c\in D$. What is the best constant
$\rho $ in the inequality%
\[
\limsup\limits_{i\rightarrow \infty }\{ \,\,\vert H_{i
}( f,c) \vert ^{1/2l_{\sigma ( i) }}\}
\leq \rho ,\ f\in A( D) ,
\]%
where $H_{i}( f,c) =\det ( a_{k( \alpha )
+k( \beta ) }) _{a,\beta =1}^{i}$ and $f( z)
=\sum_{i=1}^{\infty }a_{k( i) }\ ( z-c) ^{k(
i) }$?
\end{problem}

Our conjecture is that $\rho =d( c,\partial D) $ for good enough
domains. At least, due to Schiffer-Siciak result and \cite{Zint},
it is confirmed for polydiscs in $\mathbb{C}^{n}$.

\begin{acknowledgement}
We thank the anonymous referee whose suggestions improve the paper.
\end{acknowledgement}


\begin{thebibliography}{99}
\bibitem{APS} M. Andersson, M. Passare, and R. Sigurdsson, \textit{Complex
Convexity and analytic functionals I, \ }Report RH-06-95, Science Institute,
University Iceland, 1995.

\bibitem{Dn} S. Dineen, \textit{Complex Analysis in Locally Convex Spaces, }%
North Holland Publishing Company, 1981.

\bibitem{F} M. Fekete, \textit{\"{U}ber die Verteilung der Wurzen bei
gewissener algebraischen Gleichungen mit ganzzahligen Koefficienten, Math.
Z. }\textbf{17} (1923), 228-249.

\bibitem{G} G. M. Goluzin, \textit{Geometric Theory of Functions of a
Complex Variable,} Nauka, Moscow, 1966 (in Russian); English transl.:
Transl. Math. Monogr. \textbf{26, }Amer. Math. Soc., 1969.



\bibitem{Hrm} L.H\"{o}rmander, \textit{An Introduction to Complex Analysis in
Several Variables}, North Holland and American Elsevier, 1990.

\bibitem{Jed} M. Jedrzejowski, \textit{The Homogeneous Transfinite Diameter of a Compact Set in $\mathbb{C}^n$}, Annales Polonici Mathematici. 55(1991), 191-205.

\bibitem{Kis} C. O. Kiselman, \textit{Weak Lineal Convexity }Banach Center Publications, \textbf{107}(2016), 159-174.


\bibitem{Lj} F. Leja, \textit{Probl\`{e}mes \`{a} r\'{e}sondre pos\'{e}s
\`{a} la conf\'{e}rence, }Colloq. Math. \textbf{7 }(1959), 151-153.

\bibitem{Lv} N. Levenberg, \textit{Capacities in Several Complex Variables},
Thesis, University of Michigan, 1984.


\bibitem{Mart} A. Martineau, \textit{Sur la notion d'ensemble fortement lineelement convexe}, Anais Acado
Brasill. Cienc., \textbf{4} No. 4, 427-435 (1968).

\bibitem{P} G. Polya, \textit{Beitrag zur Verallgemeinerung des
Verzerrungssatzes auf mehrfach zusammenh\"{a}ngende Gebiete. III}, Sitzber.
Preuss. Akad. Wiss., Phys.-Math., (1929), 55--62.

\bibitem{SS} M. Schiffer, J. Siciak, \textit{Transfinite diameter and
analytic continuation of functions of two complex variables}, in: Studies
Math. Analysis and Related topics, Stanford (1962), pp.341--358.

\bibitem{Sh1} V. P. Sheinov, \textit{Transfinite diameter and some theorems
of Polya in the case of several complex variables, }Siberian Math. J.
\textbf{6 }(1971), 999-1004.

\bibitem{Sh2} V. P. Sheinov, \textit{Invariant form of Polya's inequalities}%
, Siberian Math. J. \textbf{14} (1973), 194--203 (English Translation).

\bibitem{St} E. J. Straube, \textit{Power series with integer coefficients
in several variables, C}omment Math. Helv. \textbf{62 }(1987), 602-615.


\bibitem{Z1} V. Zakharyuta, \textit{Transfinite diameter, Chebyshev
constants and capacities in }$\mathbb{C}^{n}$, Math. USSR Sbornik, \textbf{%
25 }(1975), 350--364.

\bibitem{Zint} V. Zakharyuta, \textit{Internal characteristics of domains in $\mathbb{C}^{n}$}, Annales Polonici Mathematici, (2014),Vol. 111, No.3, 215--236.

\bibitem{Zn} S. V. Znamenskii, \textit{Strong linear convexity. I. Duality
of spaces of holomorphic functions. }Siberian Math. J. \textbf{26 }(1985),
331-341.

\bibitem{Znv1} V. A. Znamenskii, \textit{The stability of the transfinite
diameter for a compact set in }$\mathbb{C}^{n}$, Izv. Severo-Kavkaz. Nauchn.
Tsentra Vyssh. Shkoly Ser. Estestv. Nauk, 1978, no. 4, 14-16.

\bibitem{Znv2} V. A. Znamenskii, \textit{Chebyshev constant for a compact
set in }$\mathbb{C}^{n}$, Turkish J. Math. \textbf{18 } no. 3\textit{\ }%
(1994), 229-237.
\end{thebibliography}
\end{document}